\documentclass[a4paper,11pt]{amsart}
\usepackage{amssymb}
\usepackage[arrow,matrix]{xy}

\title{Universal derived equivalences of posets}
\author{Sefi Ladkani}

\address{Einstein Institute of Mathematics, The Hebrew University of Jerusalem, Jerusalem 91904, Israel}
\email{sefil@math.huji.ac.il}

\DeclareMathOperator{\Hom}{Hom}
\DeclareMathOperator{\img}{im}
\DeclareMathOperator{\cone}{C}

\DeclareMathOperator{\id}{id}
\DeclareMathOperator{\hh}{H}
\DeclareMathOperator{\Func}{Func}

\newcommand{\eps}{\varepsilon}
\newcommand{\vphi}{\varphi}
\newcommand{\one}{\mathbf{1}}

\newcommand{\bZ}{\mathbb{Z}}
\newcommand{\cA}{\mathcal{A}}

\newcommand{\cB}{\mathcal{B}}
\newcommand{\cC}{\mathcal{C}}
\newcommand{\cF}{\mathcal{F}}
\newcommand{\cI}{\mathcal{I}}

\theoremstyle{plain}
\newtheorem{theorem}{Theorem}[section]
\newtheorem{lemma}[theorem]{Lemma}
\newtheorem{prop}[theorem]{Proposition}
\newtheorem{cor}[theorem]{Corollary}

\theoremstyle{definition}
\newtheorem{defn}[theorem]{Definition}

\newtheorem{example}[theorem]{Example}

\numberwithin{equation}{section}

\begin{document}

\begin{abstract}
By using only combinatorial data on two posets $X$ and $Y$, we
construct a set of so-called formulas. A formula produces
simultaneously, for any abelian category $\cA$, a functor between the
categories of complexes of diagrams over $X$ and $Y$ with values in
$\cA$. This functor induces a triangulated functor between the
corresponding derived categories.

This allows us to prove, for pairs $X, Y$ of posets sharing certain
common underlying combinatorial structure, that for any abelian
category $\cA$, regardless of its nature, the categories of diagrams
over $X$ and $Y$ with values in $\cA$ are derived equivalent.
\end{abstract}

\maketitle

\section{Introduction}

In previous work~\cite{Ladkani06} we considered the question when the
categories $\cA^X$ and $\cA^Y$ of diagrams over finite posets $X$ and
$Y$ with values in the abelian category $\cA$ of finite dimensional
vector spaces over a fixed field $k$, are derived equivalent.

Since in that case the category of diagrams $\cA^X$ is equivalent to
the category of finitely generated modules over the incidence algebra
$kX$, methods from the theory of derived equivalence of algebras, in
particular tilting theory, could be
used~\cite{Happel88,Rickard89,Rickard91}.

Interestingly, in all cases considered, the derived equivalence of two
categories of diagrams does not depend on the field $k$. A natural
question arises whether there is a general principle which explains
this fact and extends to any arbitrary abelian category $\cA$.

In this paper we provide a positive answer in the following sense; we
exhibit several constructions of pairs of posets $X$ and $Y$ such that
the derived categories $D(\cA^X)$ and $D(\cA^Y)$ are equivalent for any
abelian category $\cA$, regardless of its nature. Such pairs of posets
are called \emph{universally derived equivalent}, since the derived
equivalence is universal and originates from the combinatorial and
topological properties of the posets, rather than the specific abelian
categories involved.

Our main tools are the so-called formulas. A formula consists of
combinatorial data that produces simultaneously, for any abelian
category $\cA$, a functor between the categories of complexes of
diagrams over $X$ and $Y$ with values in $\cA$, which induces a
triangulated functor between the corresponding derived categories.

\subsection{The main construction}
Let $X$ and $Y$ be two finite partially ordered sets (\emph{posets}).
For $y \in Y$, write $[y,\cdot] = \{ y' \in Y \,:\, y' \geq y \}$ and
$[\cdot,y] = \{ y' \in Y \,:\, y' \leq y \}$. Let $\{Y_x\}_{x \in X}$
be a collection of subsets of $Y$ indexed by the elements of $X$, such
that
\begin{equation} \label{e:Yantichain}
[y,\cdot] \cap [y',\cdot] = \phi \quad \text{and} \quad [\cdot,y] \cap
[\cdot,y'] = \phi
\end{equation}
for any $x \in X$ and $y \neq y'$ in $Y_x$. Assume in addition that for
any $x \leq x'$, there exists an isomorphism $\vphi_{x,x'} : Y_x
\xrightarrow{\sim} Y_{x'}$ such that
\begin{align} \label{e:Yleq}
y &\leq \vphi_{x,x'}(y) && \text{for all $y \in Y_x$}
\end{align}
By~\eqref{e:Yantichain}, it follows that
\begin{align}\label{e:Ysystem}
\vphi_{x,x''} &= \vphi_{x',x''} \vphi_{x,x'} && \text{for all $x \leq
x' \leq x''$.}
\end{align}

Define two partial orders $\leq_{+}$ and $\leq_{-}$ on the disjoint
union $X \sqcup Y$ as follows. Inside $X$ and $Y$, the orders
$\leq_{+}$ and $\leq_{-}$ agree with the original ones, and for $x \in
X$ and $y \in Y$ we set
\begin{align} \label{e:XYorder}
x \leq_{+} y \Longleftrightarrow \text{$\exists\, y_x \in Y_x$ with
$y_x \leq y$} \\ \notag
y \leq_{-} x \Longleftrightarrow \text{$\exists\, y_x \in Y_x$ with
$y \leq y_x$}
\end{align}
with no other relations (note that the element $y_x$ is unique
by~\eqref{e:Yantichain}, and that $\leq_{+}$, $\leq_{-}$ are partial
orders by~\eqref{e:Yleq}).

\begin{theorem} \label{t:XYorder}
The two posets $(X \sqcup Y, \leq_{+})$ and $(X \sqcup Y, \leq_{-})$
are universally derived equivalent.
\end{theorem}

The assumption~\eqref{e:Yantichain} of the Theorem cannot be dropped,
as demonstrated by the following example.

\begin{example}
Consider the two posets whose Hasse diagrams are given by
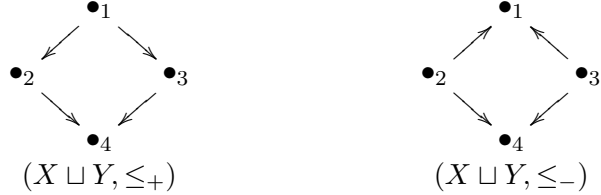
\[
\begin{array}{ccc}
\xymatrix@=1pc{
& {\bullet_1} \ar[dl] \ar[dr] \\
{\bullet_2} \ar[dr] & & {\bullet_3} \ar[dl] \\
& {\bullet_4}
}
& \hspace{60pt} &
\xymatrix@=1pc{
& {\bullet_1} \\
{\bullet_2} \ar[dr] \ar[ur] & & {\bullet_3} \ar[dl] \ar[ul] \\
& {\bullet_4} }
\\
(X \sqcup Y, \leq_{+}) & & (X \sqcup Y, \leq_{-})
\end{array}
\]
They can be represented as $(X \sqcup Y, \leq_{+})$ and $(X \sqcup Y,
\leq_{-})$ where $X = \{1\}$, $Y = \{2, 3, 4\}$ and $Y_1 = \{2, 3\}
\subset Y$. The categories of diagrams over these two posets are in
general not derived equivalent, even for diagrams of vector spaces.
\end{example}

The construction of Theorem~\ref{t:XYorder} has many interesting
consequences, some of them related to ordinal sums and others to
generalized BGP reflections~\cite{BGP73}. First, consider the case
where all the subsets $Y_x$ are single points, that is, there exists a
function $f : X \to Y$ with $Y_x = \{f(x)\}$ for all $x \in X$.
Then~\eqref{e:Yantichain} and~\eqref{e:Ysystem} are automatically
satisfied and the condition~\eqref{e:Yleq} is equivalent to $f$ being
\emph{order preserving}, i.e. $f(x) \leq f(x')$ for $x \leq x'$. Let
$\leq_{+}^f$ and $\leq_{-}^f$ denote the corresponding orders on $X
\sqcup Y$, and note that~\eqref{e:XYorder} takes the simplified form
\begin{align} \label{e:XYorderf}
x \leq_{+}^{f} y \Longleftrightarrow f(x) \leq y \\
\notag y \leq_{-}^{f} x \Longleftrightarrow y \leq f(x)
\end{align}

\begin{cor} \label{c:XYorderf}
Let $f:X \to Y$ be order preserving. Then the two posets $(X \sqcup Y,
\leq_{+}^f)$ and $(X \sqcup Y, \leq_{-}^f)$ are universally derived
equivalent.
\end{cor}

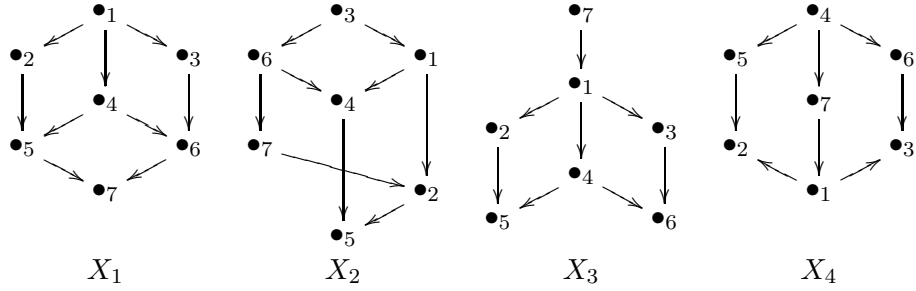
\begin{figure}
\[
\begin{array}{cccc}
\xymatrix@=0.33pc{
&& {\bullet_1} \ar[dll] \ar[dd] \ar[drr] \\
{\bullet_2} \ar[dd] && && {\bullet_3} \ar[dd] \\
&& {\bullet_4} \ar[dll] \ar[drr] \\
{\bullet_5} \ar[drr] && && {\bullet_6} \ar[dll] \\
&& {\bullet_7}
}
&
\xymatrix@=0.33pc{
&& {\bullet_3} \ar[dll] \ar[drr] \\
{\bullet_6} \ar[dd] \ar[drr] && && {\bullet_1} \ar[dll] \ar[ddd] \\
&& {\bullet_4} \ar[ddd] \\
{\bullet_7} \ar[drrrr] \\
&& && {\bullet_2} \ar[dll] \\
&& {\bullet_5}
}
&
\xymatrix@=0.33pc{
&& {\bullet_7} \ar[dd] \\ \\
&& {\bullet_1} \ar[dll] \ar[dd] \ar[drr] \\
{\bullet_2} \ar[dd] && && {\bullet_3} \ar[dd] \\
&& {\bullet_4} \ar[dll] \ar[drr] \\
{\bullet_5} && && {\bullet_6}
}
&
\xymatrix@=0.33pc{
&& {\bullet_4} \ar[dll] \ar[dd] \ar[drr] \\
{\bullet_5} \ar[dd] && && {\bullet_6} \ar[dd] \\
&& {\bullet_7} \ar[dd]  \\
{\bullet_2} && && {\bullet_3} \\
&& {\bullet_1} \ar[ull] \ar[urr]
}
\\
X_1 & X_2 & X_3 & X_4
\end{array}
\]
\caption{Four universally derived equivalent posets}
\label{f:X1234}
\end{figure}

\begin{example}
Consider the four posets $X_1$, $X_2$, $X_3$, $X_4$ whose Hasse
diagrams are drawn in Figure~\ref{f:X1234}. For any of the pairs
$(i,j)$ where $(i,j)=(1,2)$, $(1,3)$ or $(3,4)$ we find posets $X_{ij}$
and $X_{ji}$ and an order-preserving function $f_{ij} : X_{ij} \to
X_{ji}$ such that
\begin{align*}
X_i \simeq (X_{ij} \sqcup X_{ji}, \leq_{+}^{f_{ij}}) && X_j \simeq
(X_{ij} \sqcup X_{ji}, \leq_{-}^{f_{ij}})
\end{align*}
hence $X_i$ and $X_j$ are universally derived equivalent. Indeed, let
\begin{align*}
&X_{12} = \{1, 2, 4, 5 \} && X_{21} = \{3, 6, 7 \} \\
&X_{13} = \{1, 2, 3, 4, 5, 6 \} && X_{31} = \{7 \} \\
&X_{34} = \{1, 2, 3, 7 \} && X_{43} = \{4, 5, 6\}
\end{align*}
and define $f_{12} : X_{12} \to X_{21}$, $f_{13} : X_{13} \to X_{31}$
and $f_{34} : X_{34} \to X_{43}$ by
\begin{align*}
&f_{12}(1) = 3 && f_{12}(2) = f_{12}(5) = 7 && f_{12}(4) = 6 \\
&f_{13}(1) = \dots = f_{13}(6) = 7 \\
&f_{34}(1) = f_{34}(7) = 4 && f_{34}(2) = 5 && f_{34}(3) = 6
\end{align*}
\end{example}

\subsection{Applications to ordinal sums}
Recall that the \emph{ordinal sum} of two posets $(P,\leq_P)$ and
$(Q,\leq_Q)$, denoted $P \oplus Q$, is the poset $(P \sqcup Q, \leq)$
where $x \leq y$ if $x, y \in P$ and $x \leq_P y$ or $x,y \in Q$ and $x
\leq_Q y$ or $x \in P$ and $y \in Q$. Similarly, the \emph{direct sum}
$P + Q$ is the poset $(P \sqcup Q, \leq)$ where $x \leq y$ if $x, y \in
P$ and $x \leq_P y$ or $x, y \in Q$ and $x \leq_Q y$. Note that the
direct sum is commutative (up to isomorphism) but the ordinal sum is
not. Denote by $\one$ the poset consisting of one element. Taking $Y =
\one$ in Corollary~\ref{c:XYorderf}, we get the following
\begin{cor}
For any poset $X$, the posets $X \oplus \one$ and $\one \oplus X$ are
universally derived equivalent.
\end{cor}

Note that for arbitrary two posets $X$ and $Y$, it is true that for any
field $k$, the categories of diagrams of finite dimensional $k$-vector
spaces over $X \oplus Y$ and $Y \oplus X$ are derived
equivalent~\cite[Corollary~4.15]{Ladkani06}. However the proof relies
on the notion of tilting complexes and cannot be directly extended to
arbitrary abelian categories.

In Section~\ref{ssec:X1Z} we prove the following additional consequence
of Corollary~\ref{c:XYorderf} for ordinal and direct sums.

\begin{cor} \label{c:X1Z}
For any two posets $X$ and $Z$, the posets $X \oplus \one \oplus Z$ and
$\one \oplus (X + Z)$ are universally derived equivalent. Hence the
posets $X \oplus \one \oplus Z$ and $Z \oplus \one \oplus X$ are
universally derived equivalent.
\end{cor}

The result of Corollary~\ref{c:X1Z} is no longer true when $\one$ is
replaced by an arbitrary poset, even for diagrams of vector spaces,
see~\cite[Example~4.20]{Ladkani06}.

\subsection{Generalized BGP reflections}
More consequences of Theorem~\ref{t:XYorder} are obtained by
considering the case where $X = \{*\}$ is a single point, that is,
there exists a subset $Y_0 \subseteq Y$ such that~\eqref{e:Yantichain}
holds for any $y \neq y'$ in $Y_0$. Observe that
conditions~\eqref{e:Yleq} and~\eqref{e:Ysystem} automatically hold in
this case, and the two partial orders on $Y \cup \{*\}$ corresponding
to~\eqref{e:XYorder}, denoted $\leq_{+}^{Y_0}$ and $\leq_{-}^{Y_0}$,
are obtained by extending the order on $Y$ according to
\begin{align} \label{e:Ystar}
* <_{+}^{Y_0} y \Longleftrightarrow \text{$\exists\, y_0 \in Y_0$ with
$y_0 \leq y$} \\ \notag
y <_{-}^{Y_0} * \Longleftrightarrow
\text{$\exists\, y_0 \in Y_0$ with $y \leq y_0$}
\end{align}

\begin{cor} \label{c:BGP}
Let $Y_0 \subseteq Y$ be a subset satisfying~\eqref{e:Yantichain}. Then
the posets $(Y \cup \{*\}, \leq_{+}^{Y_0})$ and $(Y \cup \{*\},
\leq_{-}^{Y_0})$ are universally derived equivalent.
\end{cor}

Note that in the Hasse diagram of $\leq_{+}^{Y_0}$, the vertex $*$ is a
source which is connected to the vertices of $Y_0$, and the Hasse
diagram of $\leq_{-}^{Y_0}$ is obtained by reverting the orientations
of the arrows from $*$, making it into a sink. Thus
Corollary~\ref{c:BGP} can be considered as a generalized BGP reflection
principle.

Viewing orientations on (finite) trees as posets by setting $x \leq y$
for two vertices $x, y$ if there exists an oriented path from $x$ to
$y$, and applying a standard combinatorial argument~\cite{BGP73}, we
recover the following corollary, already known for categories of vector
spaces over a field.

\begin{cor}
Any two orientations of a tree are universally derived equivalent.
\end{cor}

\subsection{Formulas}

By using only combinatorial data on two posets $X$ and $Y$, we
construct a set of formulas $\cF_X^Y$. A formula $\pmb{\xi}$ produces
simultaneously, for any abelian category $\cA$, a functor
$F_{\pmb{\xi},\cA}$ between the categories $C(\cA^X)$ and $C(\cA^Y)$ of
complexes of diagrams over $X$ and $Y$ with values in $\cA$. This
functor induces a triangulated functor $\widetilde{F}_{\pmb{\xi},\cA}$
between the corresponding derived categories $D(\cA^X)$ and $D(\cA^Y)$
such that the following diagram is commutative
\[
\xymatrix{
C(\cA^X) \ar[r]^{F_{\pmb{\xi},\cA}} \ar[d] & C(\cA^Y) \ar[d] \\
D(\cA^X) \ar@{.>}[r]^{\widetilde{F}_{\pmb{\xi},\cA}} & D(\cA^Y) }
\]
where the vertical arrows are the canonical localizations.

We prove Theorem~\ref{t:XYorder} by exhibiting a pair of formulas
$\pmb{\xi}^{+} \in \cF_{\leq_{+}}^{\leq_{-}}$, $\pmb{\xi}^{-} \in
\cF_{\leq_{-}}^{\leq_{+}}$ and showing that for any abelian category
$\cA$, the compositions $\widetilde{F}_{\pmb{\xi}^{+},\cA}
\widetilde{F}_{\pmb{\xi}^{-},\cA}$ and
$\widetilde{F}_{\pmb{\xi}^{-},\cA} \widetilde{F}_{\pmb{\xi}^{+},\cA}$
of the corresponding triangulated functors on the derived categories
are auto-equivalences, as they are isomorphic to the translations.
Hence $\leq_{+}$ and $\leq_{-}$ are universally derived equivalent.

\section{Complexes of diagrams}

\subsection{Diagrams and sheaves}
Let $X$ be a poset and let $\cA$ be a category.

\begin{defn}
A \emph{diagram} $(A,r)$ over $X$ with values in $\cA$ consists of the
following data:
\begin{itemize}
\item
For any $x \in X$, an object $A_x$ of $\cA$
\item
For any pair $x \leq x'$, a morphism $r_{xx'}:A_x \to A_{x'}$
(\emph{restriction map})
\end{itemize}
subject to the conditions $r_{xx} = \id_{A_x}$ and
$r_{xx''}=r_{x'x''}r_{xx'}$ for all $x \leq x' \leq x''$ in $X$.

A \emph{morphism} $f:(A,r) \to (A',r')$ of diagrams consists of
morphisms $f_x : A_x \to A'_x$ for all $x \in X$, such that for any $x
\leq x'$, the diagram
\[
\xymatrix{
A_x \ar[r]^{f_x} \ar[d]_{r_{xx'}} & A'_x \ar[d]^{r'_{xx'}} \\
A_{x'} \ar[r]_{f_{x'}} & A'_{x'} }
\]
commutes.
\end{defn}

Using these definitions, we can speak of the category of diagrams over
$X$ with values in $\cA$, which will be denoted by $\cA^X$.

We can view $X$ as a small category as follows. Its objects are the
points $x \in X$, while $\Hom_X(x,x')$ is a one-element set if $x \leq
x'$ and empty otherwise. Under this viewpoint, a diagram over $X$ with
values in $\cA$ becomes a functor $A : X \to \cA$ and a morphism of
diagrams corresponds to a natural transformation, so that $\cA^X$ is
naturally identified with the category of functors $X \to \cA$. Observe
that any functor $F : \cA \to \cA'$ induces a functor $F^X : \cA^X \to
\cA'^{X}$ by the composition $F^X(A) = F \circ A$. In terms of diagrams
and morphisms, $F^X(A,r) = (FA, Fr)$ where $(FA)_x = F(A_x)$,
$(Fr)_{xx'} = F(r_{xx'})$ and $F^X(f)_x = F(f_x)$.

If $\cA$ is additive, then $\cA^X$ is additive. Assume now that $\cA$
is abelian. In this case, $\cA^X$ is also abelian, and kernels, images,
and quotients can be computed pointwise, that is, if $f : (A,r) \to
(A',r')$ is a morphism of diagrams then $(\ker f)_x = \ker f_x$, $(\img
f)_x = \img f_x$, with the restriction maps induced from $r$, $r'$. In
particular, for any $x \in X$ the evaluation functor $-_x : \cA^X \to
\cA$ taking a diagram $(A,r)$ to $A_x$ and a morphism $f=(f_x)$ to
$f_x$, is exact.

The poset $X$ admits a natural topology, whose open sets are the
subsets $U \subseteq X$ with the property that if $x \in U$ and $x \leq
x'$ then $x' \in U$. The category of diagrams over $X$ with values in
$\cA$ can then be naturally identified with the category of sheaves
over the topological space $X$ with values in $\cA$~\cite{Ladkani06}.

\subsection{Complexes and cones}

Let $\cB$ be an additive category. A \emph{complex} $(K^{\bullet},
d_K^{\bullet})$ over $\cB$ consists of objects $K^i$ for $i \in \bZ$
with morphisms $d_K^i : K^i \to K^{i+1}$ such that $d_K^{i+1} d_K^i =
0$ for all $i \in \bZ$. If $n \in \bZ$, the \emph{shift} of
$K^{\bullet}$ by $n$, denoted $K[n]^{\bullet}$, is the complex defined
by $K[n]^i = K^{i+n}$, $d_{K[n]}^i = (-1)^n d_K^{i+n}$.

Let $(K^{\bullet}, d_K^{\bullet})$, $(L^{\bullet}, d_L^{\bullet})$ be
two complexes and $f=(f^i)_{i \in \bZ}$ a collection of morphisms $f^i
: K^i \to L^i$. If $n \in \bZ$, let $f[n] = (f[n]^i)_{i \in \bZ}$ with
$f[n]^i = f^{i+n}$. Using this notation, the condition that $f$ is a
\emph{morphism} of complexes is expressed as $f[1] d_K = d_L f$.

The $\emph{cone}$ of a morphism $f : K^{\bullet} \to L^{\bullet}$,
denoted $\cone(K^{\bullet} \xrightarrow{f} L^{\bullet})$, is the
complex whose $i$-th entry equals $K^{i+1} \oplus L^i$, with the
differential
\[
d(k^{i+1}, l^i) = (-d_K^{i+1}(k^{i+1}), f^{i+1}(k^{i+1}) + d_L^i(l^i))
\]
In a more compact form, $\cone(K^{\bullet} \xrightarrow{f} L^{\bullet})
= K[1]^{\bullet} \oplus L^{\bullet}$ with the differential acting as
the matrix
\[
\begin{pmatrix}
d_K[1] & 0 \\
f[1] & d_L
\end{pmatrix}
\]
by viewing the entries as column vectors.

When $\cB$ is abelian, the \emph{$i$-th cohomology} of $(K^{\bullet},
d_K^{\bullet})$ is defined by $\hh^i(K^{\bullet}) = \ker d_K^i / \img
d_K^{i-1}$, and $(K^{\bullet}, d_K^{\bullet})$ is \emph{acyclic} if
$\hh^i(K^{\bullet})=0$ for all $i \in \bZ$. A morphism $f : K^{\bullet}
\to L^{\bullet}$ induces morphisms $\hh^i(f) : \hh^i(K^{\bullet}) \to
\hh^i(L^{\bullet})$. $f$ is called a \emph{quasi-isomorphism} if
$\hh^i(f)$ are isomorphisms for all $i \in \bZ$.

The following lemma is standard.
\begin{lemma} \label{l:Coneqis}
$f : K^{\bullet} \to L^{\bullet}$ is a quasi-isomorphism if and only if
the cone $\cone(K^{\bullet} \xrightarrow{f} L^{\bullet})$ is acyclic.
\end{lemma}

Let $C(\cB)$ denote the category of complexes over $\cB$. Denote by
$[1] : C(\cB) \to C(\cB)$ the shift functor taking a complex
$(K^{\bullet}, d_K^{\bullet})$ to $(K[1]^{\bullet},
d_{K[1]^{\bullet}})$ and a morphism $f$ to $f[1]$. Any additive functor
$G : \cB \to \cB'$ induces an additive functor $C(G) : C(\cB) \to
C(\cB')$ by sending a complex $((K^i), (d_K^i))$ to $((G(K^i)),
(G(d_K^i)))$ and a morphism $(f^i)$ to $(G(f^i))$.

\begin{lemma} \label{l:CAX}
For any additive category $\cA$ and a poset $X$, there exists an
equivalence of categories $\Phi_{X,\cA} : C(\cA^X) \simeq C(\cA)^X$
such that for any additive category $\cA'$ and an additive functor $F:
\cA \to \cA'$, the diagram
\begin{equation} \label{e:CAX}
\xymatrix{
{C(\cA)^X} \ar[r]^{\Phi_{X,\cA}}_{\sim} \ar[d]_{C(F^X)} &
{C(\cA)^X} \ar[d]^{C(F)^X} \\
{C(\cA')^X} \ar[r]^{\Phi_{X,\cA'}}_{\sim} & {C(\cA')^X}
}
\end{equation}
commutes. In other words, we can identify a complex of diagrams with a
diagram of complexes.
\end{lemma}
\begin{proof}
Let $\cA$ be additive and let $(K^{\bullet}, d^{\bullet})$ be a complex
in $C(\cA^X)$. Denote by $d^i : K^i \to K^{i+1}$ the morphisms in
$\cA^X$ and by $d^i_x : K^i_x \to K^{i+1}_x$ the morphisms on the
stalks. Let $r^i_{xy} : K^i_x \to K^i_y$ denote the restriction maps in
the diagram $K^i$.

For a morphism $f : (K^{\bullet}, d^{\bullet}) \to (L^{\bullet},
d^{\bullet})$ in $C(\cA^X)$, denote by $f^i : K^i \to L^i$ the
corresponding morphisms in $\cA^X$ and by $f^i_x : K^i_x \to L^i_x$ the
morphisms on stalks. Define a functor $\Phi: C(\cA^X) \to C(\cA)^X$ by
\begin{align*}
\Phi_{X,\cA}(K^{\bullet}, d^{\bullet}) = (\{K^{\bullet}_x\}_{x \in X},
\{r_{xy}\}) && \Phi_{X,\cA}(f) = (f_x)_{x \in X}
\end{align*}
where $(K^{\bullet}_x)^i = K^i_x$ with differential $d^{\bullet}_x =
(d^i_x)^i$, $r_{xy} = (r^i_{xy})^i : K^{\bullet}_x \to K^{\bullet}_y$
are the restriction maps, and $f_x = (f^i_x)^i : K^{\bullet}_x \to
L^{\bullet}_x$.

The commutativity of all squares in the diagram in Figure~\ref{fig:CAX}
implies that $\Phi_{X,\cA}$ is well-defined, induces the required
equivalence and that~\eqref{e:CAX} commutes.
\begin{figure}
\[
\xymatrix@=1pc{ && L^i_x \ar[rrr]^{d^i_x} \ar'[dd]^{r^i_{xy}}[ddd] &&&
L^{i+1}_x \ar[ddd]^{r^{i+1}_{xy}} \\ \\
K^i_x \ar[ddd]_{r^i_{xy}} \ar[rrr]^{d^i_x} \ar[uurr]^{f^i_x}
&&& K^{i+1}_x \ar[ddd]^{r^{i+1}_{xy}} \ar[uurr]_{f^{i+1}_x} \\
&& L^i_y \ar'[r][rrr]^{d^i_y} &&& L^{i+1}_y \\ \\
K^i_y \ar[rrr]_{d^i_y} \ar[uurr]^{f^i_y} &&& K^{i+1}_y
\ar[uurr]_{f^{i+1}_y} }
\]
\caption{} \label{fig:CAX}
\end{figure}
\end{proof}

In the sequel, $X$ is a poset, $\cA$ is an abelian category and all
complexes are in $C(\cA^X)$.

\begin{lemma} \label{l:Hx}
$\hh^i(K^{\bullet})_x = \hh^i(K^{\bullet}_x)$
\end{lemma}
\begin{proof}
Kernels and images can be computed pointwise.
\end{proof}

\begin{lemma} \label{l:Conex}
$\cone(K^{\bullet} \xrightarrow{f} L^{\bullet})_x = \cone(K^{\bullet}_x
\xrightarrow{f_x} L^{\bullet}_x)$
\end{lemma}

\begin{cor} \label{c:qisx}
Let $f : K^{\bullet} \to L^{\bullet}$ be a morphism of complexes of
diagrams. Then $f$ is a quasi-isomorphism if and only if for every $x
\in X$, $f_x : K^{\bullet}_x \to L^{\bullet}_x$ is a quasi-isomorphism.
\end{cor}
\begin{proof}
Let $x \in X$ and $i \in \bZ$. Then by Lemmas~\ref{l:Hx}
and~\ref{l:Conex},
\[
\hh^i(\cone(K^{\bullet} \xrightarrow{f} L^{\bullet}))_x =
\hh^i(\cone(K^{\bullet} \xrightarrow{f} L^{\bullet})_x) =
\hh^i(\cone(K^{\bullet}_x \xrightarrow{f_x} L^{\bullet}_x))
\]
hence $\cone(K^{\bullet} \xrightarrow{f} L^{\bullet})$ is acyclic if
and only if $\cone(K^{\bullet}_x \xrightarrow{f_x} L^{\bullet}_x)$ are
acyclic for every $x \in X$. Using Lemma~\ref{l:Coneqis}, we see that
$f$ is a quasi-isomorphism if and only if all the $f_x$ are
quasi-isomorphisms.
\end{proof}

\subsection{Universal derived equivalence}

Recall that the derived category $D(\cB)$ of an abelian category $\cB$
is obtained by formally inverting all the quasi-isomorphisms in
$C(\cB)$. It admits a structure of a triangulated category where the
distinguished triangles in $D(\cB)$ are those isomorphic to $K' \to K
\to K'' \to K'[1]$ where $0 \to K' \to K \to K'' \to 0$ is a short
exact sequence in $C(\cB)$.

\begin{defn}
Two posets $X$ and $Y$ are \emph{universally derived equivalent} if for
any abelian category $\cA$, the derived categories $D(\cA^X)$ and
$D(\cA^Y)$ are equivalent as triangulated categories.
\end{defn}

\begin{lemma}
Let $X$ and $Y$ be universally derived equivalent. Then $X^{op}$ and
$Y^{op}$ are universally derived equivalent.
\end{lemma}

\begin{lemma}
Let $X_1$, $Y_1$ and $X_2$, $Y_2$ be two pairs of universally derived
equivalent posets. Then $X_1 \times X_2$ and $Y_1 \times Y_2$ are
universally derived equivalent.
\end{lemma}

\section{Formulas}

Throughout this section, the poset $X$ is fixed.

\subsection{The category $\cC_X$}

Viewing $X \times \bZ$ as a small category with a unique map $(x,m) \to
(x',m')$ if $x \leq x'$ and $m \leq m'$ and no maps otherwise, we can
consider the additive category $\widetilde{\cC}_X$ whose objects are
finite sequences $\{(x_i,m_i)\}_{i=1}^n$ with morphisms
$\{(x_i,m_i)\}_{i=1}^n \to \{(x'_j, m'_j)\}_{j=1}^{n'}$ specified by
$n' \times n$ integer matrices $(c_{ji})_{i,j}$ satisfying $c_{ji}=0$
unless $(x_i,m_i) \leq (x'_j, m'_j)$. That is, a morphism is a formal
$\bZ$-linear combination of arrows $(x_i,m_i) \to (x'_j, m'_j)$.
Addition and composition of morphisms correspond to the usual addition
and multiplication of matrices.

To encode the fact that squares of differentials are zero, we consider
a certain quotient of $\widetilde{\cC}_X$. Namely, let
$\widetilde{\cI}_X$ be the ideal in $\widetilde{\cC}_X$ generated by
all the morphisms $(x,m) \to (x,m+2)$ for $(x,m) \in X \times \bZ$ and
let $\cC_X = \widetilde{\cC}_X / \widetilde{\cI}_X$ be the quotient.
The objects of $\cC_X$ are still sequences $\xi=\{(x_i,m_i)\}$ and the
morphisms can again be written as integer matrices, albeit not uniquely
as we ignore the entries $c_{ji}$ whenever $m'_j - m_i \geq 2$.

Define a translation functor $[1] : \cC_X \to \cC_X$ as follows. For an
object $\xi=\{(x_i,m_i)\}_{i=1}^n$, let $\xi[1] =
\{(x_i,m_i+1)\}_{i=1}^n$. For a morphism $\vphi = (c_{ji}) :
\{(x_i,m_i)\} \to \{(x'_j,m'_j)\}$, let $\vphi[1]$ be the morphism
$\{(x_i,m_i+1)\} \to \{(x'_j, m'_j+1)\}$ specified by the same matrix
$(c_{ji})$.

Let $\cA$ be an abelian category. From now on we shall denote a complex
in $C(\cA^X)$ by $K$ instead of $K^{\bullet}$, and use
Lemma~\ref{l:CAX} to identify $C(\cA^X)$ with $C(\cA)^X$. Therefore we
may think of $K$ as a diagram of complexes in $C(\cA)$ and use the
notations $K_x$, $d_x$, $r_{xx'}$ as in the proof of that lemma.

For two additive categories $\cB$ and $\cB'$, let $\Func(\cB, \cB')$
denote the category of additive functors $\cB \to \cB'$, with natural
transformations as morphisms.

\begin{prop} \label{p:etaCX}
There exists a functor $\eta : \cC_X \to \Func(C(\cA)^X, C(\cA))$
commuting with the translations.
\end{prop}
\begin{proof}
An object $\xi=\{(x_i,m_i)\}_{i=1}^n$ defines an additive functor
$F_{\xi}$ from $C(\cA)^X$ to $C(\cA)$ by sending $K \in C(\cA)^X$ and a
morphism $f : K \to K'$ to
\begin{align} \label{e:Fxi}
F_{\xi}(K) = \bigoplus_{i=1}^n K_{x_i}[m_i] & & F_{\xi}(f) =
\bigoplus_{i=1}^n f_{x_i}[m_i]
\end{align}
where the right term is the $n \times n$ diagonal matrix whose $(i,i)$
entry is $f_{x_i}[m_i] : K_{x_i}[m_i] \to K'_{x_i}[m_i]$.

To define $\eta$ on morphisms $\xi \to \xi'$, consider first the case
that $\xi=(x,m)$ and $\xi'=(x',m')$. A morphism $\vphi = (c) : (x,m)
\to (x',m')$ in $\cC_X$ is specified by an integer $c$, with $c=0$
unless $(x,m) \leq (x',m')$. Given $K \in C(\cA)^X$, define a morphism
$\eta_{\vphi}(K) : K_x[m] \to K_{x'}[m']$ by
\begin{equation} \label{e:eta1}
\eta_{\vphi}(K) =
\begin{cases}
c \cdot r_{xx'}[m] & \text{if $m'=m$ and $x' \geq x$} \\
c \cdot d_{x'}[m] r_{xx'}[m] & \text{if $m'=m+1$ and $x' \geq x$} \\
0 & \text{otherwise}
\end{cases}
\end{equation}
Then $\eta_c : F_{\xi} \to F_{\xi'}$ is a natural transformation since
the diagrams
\begin{align} \label{e:natfrd}
\xymatrix{
K_x[m] \ar[r]^{r_{xx'}[m]} \ar[d]_{f_x[m]} &
K_{x'}[m] \ar[d]^{f_{x'}[m]} \\
K'_x[m] \ar[r]_{r'_{xx'}[m]} & K'_{x'}[m]
}
& &
\xymatrix{
K_x[m] \ar[r]^{d_x[m]} \ar[d]_{f_x[m]} &
K_x[m+1] \ar[d]^{f_x[m+1]} \\
K'_x[m] \ar[r]_{d'_x[m]} & K'_x[m+1]
}
\end{align}
commute.

Let $\vphi' = (c') : (x',m') \to (x'',m'')$ be another morphism in
$\cC_X$. Then~\eqref{e:eta1} and the three relations $r_{xx''} =
r_{x'x''} r_{xx'}$, $r_{xx'}[1] d_x = d_{x'} r_{xx'}$ and $d_x[1] d_x =
0$, imply that
\begin{equation} \label{e:etafunc}
\eta_{\vphi' \vphi}(K) = \eta_{\vphi'}(K) \eta_{\vphi}(K)
\end{equation}
for every $K \in C(\cA)^X$.

Now for a general morphism $\vphi: \{(x_i,m_i)\}_{i=1}^n \to
\{(x'_j,m'_j)\}_{j=1}^{n'}$, define morphisms $\eta_{\vphi}(K) :
\bigoplus_{i=1}^{n} K_{x_i}[m_i] \to \bigoplus_{j=1}^{n'}
K_{x'_j}[m'_j]$ by
\begin{equation} \label{e:etaphi}
(\eta_{\vphi})_{ji} = \eta_{(c_{ji})} : K_{x_i}[m_i] \to K_{x'_j}[m'_j]
\end{equation}
where $\eta_{(c_{ji})}$ is defined by~\eqref{e:eta1} for $c_{ji} :
(x_i,m_i) \to (x'_j,m'_j)$.

It follows from~\eqref{e:natfrd} by linearity that for $f : K \to K'$,
\begin{equation} \label{e:etanat}
F_{\xi'}(f) \eta_{\vphi}(K) = \eta_{\vphi}(K') F_{\xi}(f)
\end{equation}
so that $\eta_{\vphi} : F_{\xi} \to F_{\xi'}$ is a natural
transformation. Linearity also shows that~\eqref{e:etafunc} holds for
general morphisms $\vphi$, $\vphi'$.

Finally, note that by~\eqref{e:Fxi} and~\eqref{e:eta1},
\begin{align*}
[1] \circ F_{\xi} = F_{\xi} \circ [1] = F_{\xi[1]} && [1] \circ
\eta_{\vphi} = \eta_{\vphi} \circ [1] = \eta_{\vphi[1]}
\end{align*}
for any object $\xi$ and morphism $\vphi$.
\end{proof}

\subsection{Formula to a point}

So far the differentials on the complexes $F_{\xi}(K)$ were just the
direct sums $\bigoplus_{i=1}^n d_{x_i}[m_i]$. For the applications,
more general differentials are needed.

Let $\vphi = (c_{ji}) : \xi \to \xi'$ be a morphism. Define
$\vphi^{\star} : \xi \to \xi'$ by $\vphi^{\star} = (c_{ji}^{\star})$
where $c_{ji}^{\star} = (-1)^{m'_j-m_i} c_{ji}$.

\begin{lemma} \label{l:etaDiff}
Let $D : \xi \to \xi[1]$ be a morphism and assume that $D^{\star}[1]
\cdot D = 0$ in $\cC_X$. Then for any $K \in C(\cA)^X$, $\eta_D(K)$ is
a differential on $F_{\xi}(K)$.
\end{lemma}
\begin{proof}
Since $F_{\xi[1]}(K) = F_{\xi}(K)[1]$, the morphism $D$ induces a map
$\eta_D(K) : F_{\xi}(K) \to F_{\xi}(K)[1]$. Thinking of $\eta_D(K)$ as
a potential differential, observe that
\begin{equation} \label{e:etaDiff}
\eta_D(K)[1] = \eta_{-D^{\star}[1]}(K)
\end{equation}

Indeed, each component $K_x[m+1] \to K_{x'}[m'+1]$ of $\eta_D(K)[1]$ is
obtained from $K_x[m] \to K_{x'}[m']$ by a change of sign. When $m'=m$,
changing the sign of a map $r_{xx'}[m]$ leads to the map
$-r_{xx'}[m+1]$. When $m'=m+1$, changing the sign of
$d_{x'}[m]r_{xx'}[m]$ leads to $d_{x'}[m+1]r_{xx'}[m+1]$, as the sign
change is already carried out in the shift of the differential
$d_{x'}[m]$. Therefore in both cases a the coefficient $c$ of $(x,m)
\to (x',m')$ changes to $-c^{\star}$.

Now the claim follows from
\[
\eta_D(K)[1] \cdot \eta_D(K) = \eta_{-D^{\star}[1]}(K) \eta_D(K) =
\eta_{- D^{\star}[1] D}(K) = 0
\]
\end{proof}

\begin{defn}
A morphism $\vphi = (c) : (x,m) \to (x',m')$ is a \emph{differential}
if $m'=m+1$, $x'=x$ and $c=1$. $\vphi$ is a \emph{restriction} if
$m'=m$ and $x' \geq x$.

A morphism $\vphi : \xi \to \xi'$ is a \emph{restriction} if all its
nonzero components are restrictions.
\end{defn}

\begin{defn}
A \emph{formula to a point} is a pair $(\xi,D)$ where
$\xi=\{(x_i,m_i)\}_{i=1}^n$ is an object of $\cC_X$ and $D =
(D_{ji})_{i,j=1}^n : \xi \to \xi[1]$ is morphism satisfying:
\begin{enumerate}
\item
$D^{\star}[1] \cdot D = 0$.

\item
$D_{ji} = 0$ for all $i > j$.

\item
$D_{ii}$ are differentials for all $1 \leq i \leq n$.
\end{enumerate}

A \emph{morphism of formulas to a point} $\vphi : (\xi,D) \to (\xi',
D')$ is a morphism $\vphi : \xi \to \xi'$ in $\cC_X$ which is a
restriction and satisfies $\vphi[1] D = D' \vphi$.
\end{defn}

Denote by $\cF_X$ the category of formulas to a point and their
morphisms. The translation $[1]$ of $\cC_X$ induces a translation $[1]$
on $\cF_X$ by $(\xi,D)[1] = (\xi[1],D[1])$ with the same action on
morphisms.

\begin{prop} \label{p:etaFX}
There exists a functor $\eta : \cF_X \to \Func(C(\cA)^X, C(\cA))$.
\end{prop}
\begin{proof}
We actually show that the required functor is induced from the functor
$\eta$ of Proposition~\ref{p:etaCX}.

An object $(\xi,D)$ defines an additive functor $F_{\xi,D} : C(\cA)^X
\to C(\cA)$ by sending $K \in C(\cA)^X$ and $f : K \to K'$ to
\begin{align*}
F_{\xi,D}(K) = F_{\xi}(K) & & F_{\xi,D}(f) = F_{\xi}(f)
\end{align*}
as in~\eqref{e:Fxi}. By Lemma~\ref{l:etaDiff}, $\eta_D(K)$ is a
differential on $F_{\xi}(K)$.

Now observe that $F_{\xi}(f)[1] \eta_D(K) = \eta_D(K') F_{\xi}(f)$
since $\eta_D : F_{\xi} \to F_{\xi[1]}$ is a natural transformation.
Therefore $F_{\xi}(f)$ is a morphism of complexes and $F_{\xi,D}$ is a
functor.

Let $\vphi : (\xi,D) \to (\xi',D')$ be a morphism in $\cF_X$. Since
$\vphi : \xi \to \xi'$ in $\cC_X$, we have a natural transformation
$\eta_{\vphi} : F_{\xi} \to F_{\xi'}$. It remains to show that
$\eta_{\vphi}(K)$ is a morphism of complexes. But the commutativity
with the differentials $\eta_D(K)$ and $\eta_{D'}(K)$ follows from
$\vphi[1] D = D' \vphi$ and the functoriality of $\eta$.
\end{proof}

\begin{example}[Zero dimensional chain] \label{ex:dim0}
Let $x \in X$ and consider $\xi = \{(x,0)\}$ with $D=(1)$. The functor
$F_{(x,0),(1)}$ sends $K$ to the stalk $K_x$ and $f : K \to K'$ to
$f_x$.
\end{example}

\begin{example}[One dimensional chain] \label{ex:dim1}
Let $x < y$ in $X$ and consider $\xi = \{(x,1),(y,0)\}$ with the map $D
= \left(\begin{smallmatrix} 1 & 0 \\ 1 & 1
\end{smallmatrix} \right) : \xi \to \xi[1]$. Then for $K \in C(\cA)^X$
and $f : K \to K'$,
\begin{align*}
F_{\xi,D}(K) = K_x[1] \oplus K_y & & F_{\xi,D}(f) = \begin{pmatrix}
f_x[1] & 0 \\ 0 & f_y
\end{pmatrix}
\end{align*}
with the differential
\[
\eta_D(K) =
\begin{pmatrix}
d_x[1] & 0 \\ r_{xy}[1] & d_y
\end{pmatrix}
: K_x[1] \oplus K_y \to K_x[2] \oplus K_y[1]
\]
Since for any object $K$, $F_{\xi,D}(K) = \cone(K_x
\xrightarrow{r_{xy}} K_y)$ as complexes, we see that for any $x < y$,
the cone $\cone(K_x \xrightarrow{r_{xy}} K_y)$ defines a functor
$C(\cA)^X \to C(\cA)$.
\end{example}

\begin{lemma} \label{l:Fxishift}
There exists a natural isomorphism $\eps : [1] \circ \eta
\xrightarrow{\simeq} \eta \circ [1]$.
\end{lemma}
\begin{proof}
We first remark that for an object $(\xi,D) \in \cF_X$, a morphism
$\vphi$ and $K \in C(\cA)^X$, $F_{\xi[1],D[1]}(K) = F_{\xi,D}(K[1])$
and $\eta_{\vphi[1]}(K) = \eta_{\vphi}(K[1])$ so that $(\eta \circ
[1])(\xi,D)$ can be viewed as first applying the shift on $C(\cA)^X$
and then applying $F_{\xi,D}$.

We will construct natural isomorphisms of functors $\eps_{\xi,D} : [1]
\circ F_{\xi,D} \xrightarrow{\simeq} F_{\xi,D} \circ [1]$ such that the
diagrams
\begin{equation} \label{e:Fxishift}
\xymatrix{ F_{\xi,D}(K)[1] \ar[d]_{[1] \circ \eta_{\vphi}}
\ar[r]^{\eps_{\xi,D}} & F_{\xi[1],D[1]}(K) \ar[d]^{\eta_{\vphi[1]}}
\\ F_{\xi',D'}(K)[1] \ar[r]_{\eps_{\xi',D'}} & F_{\xi'[1],D'[1]}(K)}
\end{equation}
commute for all $K \in C(\cA)^X$.

By~\eqref{e:etaDiff}, $[1] \circ F_{\xi,D} = F_{\xi[1],
-D^{\star}[1]}$. Write $\xi = \{(x_i,m_i)\}_{i=1}^n$,
$D=(D_{ji})_{i,j=1}^n$, and let $I_{\xi} : \xi \to \xi$ be the morphism
defined by the diagonal matrix whose $(i,i)$ entry is $(-1)^{m_i}$. By
definition, $D_{ji}^{\star} = (-1)^{m_j+1-m_i} D_{ji}$, or equivalently
$(-1)^{m_j} D_{ji} = -D^{\star}_{ji} (-1)^{m_i}$ for all $i, j$, hence
$I_{\xi}[1] D = -D^{\star} I_{\xi}$. Therefore $I_{\xi}[1] : (\xi[1],
D[1]) \to (\xi[1], -D^{\star}[1])$ is an isomorphism in $\cF_X$, so we
define $\eps_{\xi,D} = \eta_{I_{\xi}[1]}$.

For the commutativity of~\eqref{e:Fxishift}, first observe that $[1]
\circ \eta_{\vphi} = \eta_{\vphi} \circ [1] = \eta_{\vphi[1]}$. Now use
the fact that $I_{\xi'} \vphi = \vphi I_{\xi}$ for any restriction
$\vphi : \xi \to \xi'$.
\end{proof}

In the next few lemmas, we fix a formula to a point $(\xi,D)$.

\begin{lemma} \label{l:Fxises}
$F_{\xi,D}$ maps short exact sequences to short exact sequences.
\end{lemma}
\begin{proof}
Write $\xi = \{(x_i,m_i)\}_{i=1}^n$ and let $0 \to K' \xrightarrow{f'}
K \xrightarrow{f''} K'' \to 0$ be a short exact sequence. Then $0 \to
K'_x \xrightarrow{f'_x} K_x \xrightarrow{f''_x} K''_x \to 0$ is exact
for any $x \in X$, hence
\[
0 \to \bigoplus_{i=1}^n K'_{x_i}[m_i] \xrightarrow{\bigoplus_{i=1}^n
f'_{x_i}[m_i]} \bigoplus_{i=1}^n K_{x_i}[m_i]
\xrightarrow{\bigoplus_{i=1}^n f''_{x_i}[m_i]} \bigoplus_{i=1}^n
K''_{x_i}[m_i] \to 0
\]
is exact.
\end{proof}

By composing with the equivalence $\Phi : C(\cA^X) \to C(\cA)^X$, we
may view $F_{\xi,D}$ as a functor $C(\cA^X) \to C(\cA)$ between two
categories of complexes.

\begin{lemma} \label{l:Fxiqis}
$F_{\xi,D}$ maps quasi-isomorphisms to quasi-isomorphisms.
\end{lemma}
\begin{proof}
Write $\xi = \{(x_i,m_i)\}_{i=1}^n$. We prove the claim by induction on
$n$. When $n=1$, we have $\xi = (x,m)$, $F_{\xi,D}(K) = K_x[m]$ and
$F_{\xi,D}(f) = f_x[m]$, so that the claim follows from
Corollary~\ref{c:qisx}.

Assume now that $n > 1$, and let $\xi' = \{(x_i,m_i)\}_{i=1}^{n-1}$ and
$D'=(D_{ji})_{i,j=1}^{n-1}$ be the corresponding restricted matrix. By
the assumption that $D = (D_{ji})$ is lower triangular with ones on the
main diagonal, we have that the canonical embedding $\iota_K :
K_{x_n}[m_n] \to \bigoplus_{i=1}^n K_{x_i}[m_i]$ and the projection
$\pi_K : \bigoplus_{i=1}^n K_{x_i}[m_i] \to \bigoplus_{i=1}^{n-1}
K_{x_i}[m_i]$ commute with the differentials, hence there exists a
functorial short exact sequence
\begin{equation} \label{e:Fxises}
0 \to (K_{x_n}[m_n], d_{x_n}[m_n]) \to (F_{\xi,D}(K), \eta_D(K)) \to
(F_{\xi',D'}(K), \eta_{D'}(K)) \to 0
\end{equation}

Let $f: K \to K'$ be a morphism. The functoriality of~\eqref{e:Fxises}
gives rise to the following diagram of long exact sequences in
cohomology,
\[
\xymatrix{
\ar[r] & \hh^i(F_{\xi',D'}(K)) \ar[d]^{\hh^i(F_{\xi',D'}(f))} \ar[r] &
\hh^i(K_{x_n}[m_n]) \ar[d]^{\hh^i(f_{x_n}[m_n])} \ar[r] &
\hh^i(F_{\xi,D}(K))= \ar[d]^{\hh^i(F_{\xi,D}(f))} & \\
\ar[r] & \hh^i(F_{\xi',D'}(K')) \ar[r] & \hh^i(K'_{x_n}[m_n]) \ar[r]
& \hh^i(F_{\xi,D}(K'))= & \\
& =\hh^i(F_{\xi,D}(K)) \ar[d]^{\hh^i(F_{\xi,D}(f))} \ar[r] &
\hh^{i+1}(F_{\xi',D'}(K)) \ar[d]^{\hh^{i+1}(F_{\xi',D'}(f))} \ar[r] &
\hh^{i+1}(K_{x_n}[m_n]) \ar[d]^{\hh^{i+1}(f_{x_n}[m_n])} \ar[r] & \\
& =\hh^i(F_{\xi,D}(K')) \ar[r] & \hh^{i+1}(F_{\xi',D'}(K')) \ar[r] &
\hh^{i+1}(K'_{x_n}[m_n]) \ar[r] &
}
\]

Now assume that $f : K \to K'$ is a quasi-isomorphism. By the induction
hypothesis, $f_{x_n}[m_n] : K_{x_n}[m_n] \to K'_{x_n}[m_n]$ and
$F_{\xi',D'}(f) : F_{\xi',D'}(K) \to F_{\xi',D'}(K')$ are
quasi-isomorphisms, hence by the Five Lemma, $F_{\xi,D}(f)$ is also a
quasi-isomorphism.
\end{proof}

\begin{cor}
Let $(\xi,D)$ be a formula to a point. Then $F_{\xi,D}$ induces a
triangulated functor $\widetilde{F}_{\xi,D} : D(\cA^X) \to D(\cA)$.
\end{cor}

\subsection{General formulas}

\begin{defn}
Let $Y$ be a poset. A \emph{formula from $X$ to $Y$} is a diagram over
$Y$ with values in $\cF_X$.
\end{defn}

\begin{prop}
There exists a functor $\eta : \cF_X^Y \to \Func(C(\cA)^X, C(\cA)^Y)$.
\end{prop}
\begin{proof}
Let $\eta : \cF_X \to \Func(C(\cA)^X, C(\cA))$ be the functor of
Proposition~\ref{p:etaFX}. Then
\[
\eta^{Y} : \cF_X^Y \to \Func(C(\cA)^X, C(\cA))^Y \simeq \Func(C(\cA)^X,
C(\cA)^Y)
\]
is the required functor.
\end{proof}

Let $\pmb{\xi} \in \cF_X$ be a formula.

\begin{lemma}
$F_{\pmb{\xi}}$ maps short exact sequences to short exact sequences.
\end{lemma}
\begin{proof}
It is enough to consider each component of $F_{\pmb{\xi}}$ separately.
The claim now follows from Lemma~\ref{l:Fxises}.
\end{proof}

By composing from the left with the equivalence $\Phi : C(\cA^X) \to
C(\cA)^X$ and from the right with $\Phi^{-1} : C(\cA)^Y \to C(\cA^Y)$
we may view $F_{\pmb{\xi}}$ as a functor $C(\cA^X) \to C(\cA^Y)$
between two categories of complexes.

\begin{lemma}
$F_{\pmb{\xi}}$ maps quasi-isomorphisms to quasi-isomorphisms.
\end{lemma}
\begin{proof}
Let $f:K \to K'$ be a quasi-isomorphism. By Corollary~\ref{c:qisx}, it
is enough to show that each component of $F_{\pmb{\xi}}(f)$ is a
quasi-isomorphism in $C(\cA)$. But this follows from
Lemma~\ref{l:Fxiqis}.
\end{proof}

\begin{cor}
Let $\pmb{\xi}$ be a formula. Then $F_{\pmb{\xi}}$ induces a
triangulated functor $\widetilde{F}_{\pmb{\xi}} : D(\cA^X) \to
D(\cA^Y)$.
\end{cor}

\section{Applications of formulas}

\subsection{The chain with two elements}

As a first application we consider the case where the poset $X$ is a
chain of two elements
\[
\xymatrix{{\bullet_1} \ar[r] & {\bullet_2}}
\]

We focus on this simple case as the fundamental underlying principle of
Theorem~\ref{t:XYorder} can already be effectively demonstrated in that
case.

Let $(\xi_1,D_1)$, $(\xi_2, D_2)$ and $(\xi_{12}, D_{12})$ be the
following three formulas to a point in $\cF_{1 \to 2}$.
\begin{align} \label{e:xi12}
\xi_1 = (1,1), D_1 = (1) && \xi_{12} = ((1,1),(2,0)), D_{12} =
\begin{pmatrix} 1 & 0 \\ 1 & 1 \end{pmatrix} \\ \notag
\xi_2 = (2,0), D_2 = (1)
\end{align}

Let $\cA$ be an abelian category and $K = K_1 \xrightarrow{r_{12}} K_2$
be an object of $C(\cA^{1 \to 2}) \simeq C(\cA)^{1 \to 2}$. In the more
familiar notation,
\begin{align*}
F_{\xi_1,D_1}(K) = K_1[1] && F_{\xi_2,D_2}(K) = K_2 &&
F_{\xi_{12},D_{12}}(K) = \cone(K_1 \xrightarrow{r_{12}} K_2)
\end{align*}
see Examples~\ref{ex:dim0} and~\ref{ex:dim1}.

The morphisms
\begin{align*}
\vphi_1 = \begin{pmatrix} 1 & 0 \end{pmatrix}: \xi_{12} \to \xi_1 &&
\vphi_2 = \begin{pmatrix} 0 \\ 1 \end{pmatrix}: \xi_2 \to \xi_{12}
\end{align*}
are restrictions that satisfy $\vphi_1 D_{12} = D_1 \vphi_1$ and
$\vphi_2 D_2 = D_{12} \vphi_2$, hence
\begin{align*}
\pmb{\xi}^{-} = \xymatrix{{(\xi_{12}, D_{12})} \ar[r]^{\vphi_1} &
{(\xi_1, D_1)}} && \pmb{\xi}^{+} = \xymatrix{{(\xi_2, D_2)}
\ar[r]^{\vphi_2} & {(\xi_{12}, D_{12})}}
\end{align*}
are diagrams over $1 \to 2$ with values in $\cF_{1 \to 2}$, thus they
define functors $R^{-}, R^{+} : C(\cA^{1 \to 2}) \to C(\cA^{1 \to 2})$
inducing triangulated functors $\widetilde{R}^{-}, \widetilde{R}^{+} :
D(\cA^{1 \to 2}) \to D(\cA^{1 \to 2})$. Their values on objects $K \in
C(\cA^{1 \to 2})$ are
\begin{align} \label{e:R12}
R^{-}(K) &= \cone(K_1 \xrightarrow{r_{12}} K_2) \xrightarrow{\left(
\begin{smallmatrix} r_{11}[1] & 0 \end{smallmatrix} \right)} K_1[1] \\
\notag R^{+}(K) &= K_2 \xrightarrow{\left( \begin{smallmatrix} 0 \\
r_{22}
\end{smallmatrix} \right)} \cone(K_1 \xrightarrow{r_{12}} K_2)
\end{align}

\begin{prop} \label{p:Rpm1}
There are natural transformations
\[
R^{+} \circ R^{-} \xrightarrow{\eps^{+-}} [1] \xrightarrow{\eps^{-+}}
R^{-} \circ R^{+}
\]
such that $\eps^{+-}(K)$, $\eps^{-+}(K)$ are quasi-isomorphisms for all
$K \in C(\cA^{1 \to 2})$.
\end{prop}
\begin{proof}
The functors $R^{+} \circ R^{-}$ and $R^{-} \circ R^{+}$ correspond to
the compositions $\pmb{\xi}^{+-} = \pmb{\xi}^{+} \circ (\pmb{\xi}^{-}_1
\to \pmb{\xi}^{-}_2)$ and $\pmb{\xi}^{-+} = \pmb{\xi}^{-} \circ
(\pmb{\xi}^{+}_1 \to \pmb{\xi}^{+}_2)$, given by
\begin{align*}
\pmb{\xi}^{+-} &= (\xi_1, D_1)
\xrightarrow{\left(\begin{smallmatrix} 0 \\ 0 \\ 1
\end{smallmatrix}\right)}
(\xi_{121}, D_{121}) \\
\pmb{\xi}^{-+} &= (\xi_{212}, D_{212})
\xrightarrow{\left(\begin{smallmatrix} 1 & 0 & 0
\end{smallmatrix}\right)}
(\xi_2[1], D_2[1])
\end{align*}
where
\begin{align} \label{e:D121}
(\xi_{121}, D_{121}) &= \Bigl( ((1,2),(2,1),(1,1)),
\begin{pmatrix} 1 & 0 & 0 \\ -1 & 1 & 0 \\ 1 & 0 & 1 \end{pmatrix}
\Bigr) \\ \notag
(\xi_{212}, D_{212}) &= \Bigl( ((2,1),(1,1),(2,0)),
\begin{pmatrix} 1 & 0 & 0 \\ 0 & 1 & 0 \\ 1 & 1 & 1 \end{pmatrix}
\Bigr)
\end{align}
and the translation $[1]$ corresponds to the diagram
\[
\pmb{\nu} = (\xi_1, D_1) \xrightarrow{\left( \begin{smallmatrix} 1
\end{smallmatrix} \right)} (\xi_2[1], D_2[1])
\]

Let $\alpha_1, \alpha_2, \beta_1, \beta_2$ be the morphisms
\begin{align} \label{e:alphabeta}
&\alpha_1 : (\xi_1, D_1) \xrightarrow{\left(
\begin{smallmatrix} 1 \\ -1 \\ 0 \end{smallmatrix} \right)}
(\xi_{212}, D_{212})
&& \beta_1 :(\xi_{212}, D_{212}) \xrightarrow{\left(
\begin{smallmatrix} 0 & -1 & 0 \end{smallmatrix} \right)}
(\xi_1, D_1) \\ \notag
&\alpha_2 : (\xi_2[1], D_2[1]) \xrightarrow{\left(
\begin{smallmatrix} 0 \\ 1 \\ 0 \end{smallmatrix} \right)}
(\xi_{121}, D_{121})
&& \beta_2 : (\xi_{121}, D_{121}) \xrightarrow{\left(
\begin{smallmatrix} 0 & 1 & 1 \end{smallmatrix} \right)}
(\xi_2[1], D_2[1])
\end{align}
The following diagram in $\cF_{1 \to 2}$
\[
\xymatrix{
(\xi_1, D_1)
\ar[r]^{\left( \begin{smallmatrix} 1 \end{smallmatrix} \right)}
\ar[d]_{\left( \begin{smallmatrix} 0 \\ 0 \\ 1 \end{smallmatrix} \right)}
& (\xi_1, D_1)
\ar[r]^{\alpha_1}
\ar[d]^{\left( \begin{smallmatrix} 1 \end{smallmatrix} \right)}
& (\xi_{212}, D_{212})
\ar[d]^{\left( \begin{smallmatrix} 1 & 0 & 0 \end{smallmatrix} \right)}
\\
(\xi_{121}, D_{121})
\ar[r]_{\beta_2}
& (\xi_2[1], D_2[1])
\ar[r]^{\left( \begin{smallmatrix} 1 \end{smallmatrix} \right)}
& (\xi_2[1], D_2[1])
}
\]
is commutative, hence the horizontal arrows induce morphisms of
formulas $\pmb{\xi}^{+-} \to \pmb{\nu}$ and $\pmb{\nu} \to
\pmb{\xi}^{-+}$, inducing natural transformations $\eps^{+-} : R^{+}
R^{-} \to [1]$ and $\eps^{-+} : [1] \to R^{-} R^{+}$.

We prove that $\eps^{+-}(K)$ and $\eps^{-+}(K)$ are quasi-isomorphisms
for all $K$ by showing that each component is a quasi-isomorphism (see
Corollary~\ref{c:qisx}). Indeed, let $h_1 : \xi_{212} \to
\xi_{212}[-1]$ and $h_2 : \xi_{121} \to \xi_{121}[-1]$ be the maps
\[
h_1 = h_2 = \begin{pmatrix} 0 & 0 & 1 \\ 0 & 0 & 0 \\ 0 & 0 & 0
\end{pmatrix}
\]
Then
\begin{align} \label{e:abhomotopy}
\beta_1 \alpha_1 = \begin{pmatrix} 1 \end{pmatrix} && \alpha_1 \beta_1
+ (h_1[1]D_{212} + D^{\star}_{212}[-1] h_1) = I_3 \\ \notag
\beta_2 \alpha_2 = \begin{pmatrix} 1 \end{pmatrix} && \alpha_2 \beta_2
+ (h_2[1]D_{121} + D^{\star}_{121}[-1] h_2) = I_3
\end{align}
where $I_3$ is the $3 \times 3$ identity matrix, hence $\beta_1
\alpha_1$ and $\beta_2 \alpha_2$ induce the identities and $\alpha_1
\beta_1$, $\alpha_2 \beta_2$ induce morphisms $\eta_{\alpha_1
\beta_1}(K)$ and $\eta_{\alpha_2 \beta_2}(K)$ homotopic to the
identities. Therefore $\eta_{\alpha_1}(K)$, $\eta_{\alpha_2}(K)$,
$\eta_{\beta_1}(K)$ and $\eta_{\beta_2}(K)$ are quasi-isomorphisms.
\end{proof}

\begin{prop}
There are natural transformations
\begin{align*}
R^{+} \circ R^{+} \xrightarrow{\eps^{++}} R^{-} && R^{+} \circ [1]
\xrightarrow{\eps^{--}} R^{-} \circ R^{-}
\end{align*}
such that $\eps^{++}(K)$, $\eps^{--}(K)$ are quasi-isomorphisms for all
$K \in C(\cA^{1 \to 2})$.
\end{prop}
\begin{proof}
The functors $R^{+} \circ R^{+}$ and $R^{-} \circ R^{-}$ correspond to
the compositions $\pmb{\xi}^{++} = \pmb{\xi}^{+} \circ (\pmb{\xi}^{+}_1
\to \pmb{\xi}^{+}_2)$ and $\pmb{\xi}^{--} = \pmb{\xi}^{-} \circ
(\pmb{\xi}^{-}_1 \to \pmb{\xi}^{-}_2)$, given by
\begin{align*}
\pmb{\xi}^{++} &= (\xi_{12}, D_{12})
\xrightarrow{\left(
\begin{smallmatrix} 0 & 0 \\ 1 & 0 \\ 0 & 1 \end{smallmatrix}
\right)} (\xi_{212}, D_{212}) \\
\pmb{\xi}^{--} &= (\xi_{121}, D_{121})
\xrightarrow{\left(
\begin{smallmatrix} 1 & 0 & 0 \\ 0 & 1 & 0 \end{smallmatrix}
\right)} (\xi_{12}[1], -D^{\star}_{12}[1])
\end{align*}
where $(\xi_{121}, D_{121})$ and $(\xi_{212}, D_{212})$ are as
in~\eqref{e:D121}. The commutative diagrams
\begin{align*}
\xymatrix{
(\xi_{12}, D_{12})
\ar[r]^{\left( \begin{smallmatrix} 1 & 0 \\ 0 & 1
\end{smallmatrix}\right)}
\ar[d]_{\left( \begin{smallmatrix} 0 & 0 \\ 1 & 0 \\ 0 & 1
\end{smallmatrix} \right)}
& (\xi_{12}, D_{12})
\ar[d]^{\left( \begin{smallmatrix} 1 & 0 \end{smallmatrix} \right)}
\\
(\xi_{212}, D_{212}) \ar[r]^{-\beta_1}
& (\xi_1[1], D_1[1])
}
& &
\xymatrix{
(\xi_2[1], D_2[1])
\ar[r]^{\alpha_2}
\ar[d]_{\left( \begin{smallmatrix} 0 \\ 1 \end{smallmatrix} \right)}
& (\xi_{121}, D_{121})
\ar[d]^{\left( \begin{smallmatrix} 1 & 0 & 0 \\ 0 & 1 & 0
\end{smallmatrix} \right)}
\\
(\xi_{12}[1], D_{12}[1])
\ar[r]^{\left( \begin{smallmatrix} -1 & 0 \\ 0 & 1 \end{smallmatrix}
\right)}
& (\xi_{12}[1], -D^{\star}_{12}[1])
}
\end{align*}
where $\alpha_2, \beta_1$ are as in~\eqref{e:alphabeta}, define
morphisms of formulas $\pmb{\xi}^{++} \to \pmb{\xi}^{-}$ and
$\pmb{\xi}^{+}[1] \to \pmb{\xi}^{--}$, hence natural transformations
$\eps^{++} : R^{+} R^{+} \to R^{-}$ and $\eps^{--} : R^{+}[1] \to R^{-}
R^{-}$. Using the homotopies~\eqref{e:abhomotopy}, one proves that
$\eps^{++}(K)$ and $\eps^{--}(K)$ are quasi-isomorphisms for all $K$ in
the same way as before.
\end{proof}

\begin{cor}
For any abelian category $\cA$, the functors $\widetilde{R}^{+}$ and
$\widetilde{R}^{-}$ are auto-equivalences of $D(\cA)^{1 \to 2}$
satisfying
\begin{align*}
\widetilde{R}^{+} \widetilde{R}^{-} \simeq [1] \simeq \widetilde{R}^{-}
\widetilde{R}^{+} &&
(\widetilde{R}^{+})^2 \simeq \widetilde{R}^{-} &&
(\widetilde{R}^{-})^2 \simeq \widetilde{R}^{+} \circ [1]
\end{align*}
hence $(\widetilde{R}^{+})^3 \simeq [1]$.
\end{cor}

\subsection{Proof of Theorem~\ref{t:XYorder}}

Let $X$ and $Y$ be two posets satisfying the
assumptions~\eqref{e:Yantichain} and~\eqref{e:Yleq}, and let
$\leq_{+}$, $\leq_{-}$ be the partial orders on $X \sqcup Y$ as defined
by~\eqref{e:XYorder}. We will prove the universal derived equivalence
of $\leq_{+}$ and $\leq_{-}$ by defining two formulas $\pmb{\xi}^+$,
$\pmb{\xi}^-$ that will induce, for any abelian category $\cA$,
functors
\begin{align*}
R^{+} = F_{\pmb{\xi}^{+}} : C(\cA)^{\leq_{+}} \to C(\cA)^{\leq_{-}} &&
R^{-} = F_{\pmb{\xi}^{-}} : C(\cA)^{\leq_{-}} \to C(\cA)^{\leq_{+}} \\
\intertext{and} \widetilde{R}^{+} = \widetilde{F}_{\pmb{\xi}^{+}} :
D(\cA^{\leq_{+}}) \to D(\cA^{\leq_{-}}) && \widetilde{R}^{-} =
\widetilde{F}_{\pmb{\xi}^{-}} : D(\cA^{\leq_{-}}) \to D(\cA^{\leq_{+}})
\end{align*}
such that $\widetilde{R}^{+} \widetilde{R}^{-} \simeq [1]$ and
$\widetilde{R}^{-} \widetilde{R}^{+} \simeq [1]$.

\subsubsection{Definition of the formulas to points}
For $x \in X$ and $y \in Y$, let
\begin{align*}
\xi_x = \left((x,0), \begin{pmatrix} 1 \end{pmatrix} \right) && \xi_y =
\left((y,0), \begin{pmatrix} 1 \end{pmatrix} \right) && \xi_{Y_x} =
\left((y,0)_{y \in Y_x}, I \right)
\end{align*}
where $I$ is the identity matrix. We consider $\xi_x$, $\xi_y$ and
$\xi_{Y_x}$ as formulas either in $\cF_{\leq_{+}}$ or in
$\cF_{\leq_{-}}$, as appropriate. If $y \in Y$, define
\begin{align*}
\pmb{\xi}^{+}_y = \xi_y \in \cF_{\leq_{+}} && \pmb{\xi}^{-}_y = \xi_y
\in \cF_{\leq_{-}}
\end{align*}
as in Example~\ref{ex:dim0}. If $x \in X$, let
\begin{align*}
\xi_{x, Y_x} = \bigl( \xi_x \xrightarrow{\left(
\begin{smallmatrix} 1 \\ 1 \\ \dots \\ 1 \end{smallmatrix} \right)}
\xi_{Y_x} \bigr) \in \cF_{\leq_{+}}^{1 \to 2} && \xi_{Y_x, x} = \bigl(
\xi_{Y_x} \xrightarrow{\left(
\begin{smallmatrix} 1 & 1 & \dots & 1
\end{smallmatrix}\right)} \xi_x \bigr) \in \cF_{\leq_{-}}^{1 \to 2}
\end{align*}
be formulas to $1 \to 2$ and define
\begin{align*}
\pmb{\xi}^{+}_x = \xi_{12} \circ \xi_{x, Y_x} && \pmb{\xi}^{-}_x =
\xi_{12} \circ \xi_{Y_x, x}
\end{align*}
as compositions with the formula $\xi_{12}$ defined in~\eqref{e:xi12}.

In explicit terms, let $K \in C(\cA)^{\leq_{+}}$, $L \in
C(\cA)^{\leq_{-}}$, and denote by $\{r_{xy}\}$ the restriction maps in
$K$ and by $\{s_{yx}\}$ the restriction maps in $L$. For $x \in X$ and
$y \in Y_x$, let $\iota_y : K_y \to \bigoplus_{y_x \in Y_x} K_{y_x}$
and $\pi_y : \bigoplus_{y_x \in Y_x} L_{y_x} \to L_y$ be the canonical
inclusions and projections. Then
\begin{align*}
&R^{+}(K)_x = \cone(K_x \xrightarrow{\sum_{y \in Y_x} \iota_y r_{xy}}
\bigoplus_{y \in Y_x} K_y) && R^{+}(K)_y = K_y \\
&R^{-}(L)_x = \cone(\bigoplus_{y \in Y_x} L_y \xrightarrow{\sum_{y \in
Y_x} s_{yx} \pi_y} L_x) && R^{-}(L)_y = L_y[1]
\end{align*}
for $x \in X$, $y \in Y$.

\subsubsection{Definition of the restriction maps}
We shall denote by $\rho^{+}$ the restriction maps between the formulas
in $R^{+}$ and by $\rho^{-}$ the maps between those in $R^{-}$. We
consider several cases, and use the explicit notation.

For $y \leq y'$, define
\begin{align*}
\rho^{+}_{yy'}(K) = r_{yy'} : K_y \to K_{y'} && \rho^{-}_{yy'}(L) =
s_{yy'}[1] : L_y[1] \to L_{y'}[1]
\end{align*}

For $x \leq x'$, we use the isomorphism $\vphi_{x,x'} : Y_x \to Y_{x'}$
and the property that $y \leq \vphi_{x,x'}(y)$ for all $y \in Y_x$ to
define the diagonal maps
\begin{align*}
\rho^{+}_{xx'}(K) &= r_{xx'}[1] \oplus (\bigoplus_{y \in Y_x}
r_{y,\vphi_{xx'}(y)}) : R^{+}(K)_x \to R^{+}(K)_{x'}\\
\rho^{-}_{xx'}(L) &= (\bigoplus_{y \in Y_x} s_{y,\vphi_{xx'}(y)}[1])
\oplus s_{xx'} : R^{-}(L)_x \to R^{-}(L)_{x'}
\end{align*}

If $y_x \in Y_x$, then by~\eqref{e:XYorder}, $y_x \leq_{-} x$, $x
\leq_{+} y_x$, and we define
\begin{align*}
\rho^{+}_{y_x x}(K) &= K_{y_x} \xrightarrow{\left(
\begin{smallmatrix} 0 \\ \iota_{y_x} \end{smallmatrix} \right)}
\cone (K_x \to \bigoplus_{y \in Y_x} K_{y}) \\
\rho^{-}_{x y_x}(L) &= \cone(\bigoplus_{y \in Y_x} L_y \to K_x)
\xrightarrow{\left( \begin{smallmatrix} \pi_{y_x}[1] & 0
\end{smallmatrix}\right)} L_{y_x}[1]
\end{align*}

Finally, if $y \leq_{-} x$, by~\eqref{e:Yantichain} there exists a
\emph{unique} $y_x \in Y_x$ such that $y \leq y_x$ and we set
$\rho^{+}_{y x}(K) = \rho^{+}_{y_x x}(K) \rho^{+}_{y y_x}(K)$.
Similarly, if $x \leq_{+} y$, there exists a unique $y_x \in Y_x$ with
$y_x \leq y$, and we set $\rho^{-}_{xy}(L) = \rho^{-}_{y_x y}(L)
\rho^{-}_{x y_x}(L)$.

\subsubsection{Verification of commutativity}
Again there are several cases to consider. First, when $y \leq y' \leq
y''$, $\rho^{+}_{yy''} = \rho^{+}_{y'y''} \rho^{+}_{yy'}$ follows from
the commutativity of the restrictions $r_{yy''} = r_{y'y''} r_{yy'}$,
and similarly for $\rho^{-}$.

Let $x \leq x' \leq x''$. Since $\vphi_{xx'} : Y_x \to Y_{x'}$ is an
isomorphism and $\vphi_{xx''} = \vphi_{x'x''} \vphi_{xx'}$, we can
write
\begin{align*}
\rho^{+}_{x'x''}(K) &= r_{x'x''}[1] \oplus \bigoplus_{y' \in Y_{x'}}
r_{y',\vphi_{x'x''}(y')} =  r_{x'x''}[1] \oplus \bigoplus_{y \in Y_x}
r_{\vphi_{xx'}(y), \vphi_{x'x''} \vphi_{xx'}(y)} \\
&= r_{x'x''}[1] \oplus \bigoplus_{y \in Y_x} r_{\vphi_{xx'}(y),
\vphi_{xx''}(y)}
\end{align*}
Now $\rho^{+}_{xx''} = \rho^{+}_{x'x''} \rho^{+}_{xx'}$ follows from
the commutativity of the restrictions $r_{xx''} = r_{x'x''} r_{xx'}$
and $r_{y, \vphi_{xx''}(y)} = r_{\vphi_{xx'}(y), \vphi_{xx''}(y)} r_{y,
\vphi_{xx'}(y)}$. The proof for $\rho^{-}$ is similar.

If $y' \leq y \leq_{-} x$, let $y_x, y'_x \in Y_x$ be the elements
satisfying $y \leq y_x$, $y ' \leq y'_x$. Then $y'_x = y_x$ by
uniqueness, since $y' \leq y_x$. Hence
\[
\rho^{+}_{y'x} = \rho^{+}_{y_x x} \rho^{+}_{y'y_x} = \rho^{+}_{y_x x}
\rho^{+}_{yy_x} \rho^{+}_{y'y} = \rho^{+}_{yx} \rho^{+}_{y'y}
\]
The proof for $\rho^{-}$ in the case $x \leq_{+} y \leq y'$ is similar.

If $y_x \leq_{-} x \leq x'$ where $y_x \in Y_x$, then $y_{x'} =
\vphi_{xx'}(y_x)$ is the unique element $y_{x'} \in Y_{x'}$ with $y_x
\leq y_{x'}$, and
\[
\rho^{+}_{y_x x'} = \rho^{+}_{\vphi_{xx'}(y_x), x'} \rho^{+}_{y_x,
\vphi_{xx'}(y_x)} = \rho^{+}_{xx'} \rho^{+}_{y_x x}
\]
by the commutativity of the diagram
\[
\xymatrix{
K_{y_x} \ar[rr]^(0.4){\rho^{+}_{y_x,x}}
\ar[d]_{r_{y_x,\vphi_{xx'}(y_x)}}
&& \cone(K_x \to \bigoplus_{y \in Y_x} K_y)
\ar[d]^{\rho^{+}_{xx'} =
r_{xx'}[1] \oplus \bigoplus r_{y,\vphi_{xx'}(y)}}
\\
K_{\vphi_{x,x'}(y_x)} \ar[rr]^(0.4){\rho^{+}_{\vphi_{xx'}(y_x),x'}}
&& \cone(K_{x'} \to \bigoplus_{y' \in Y_{x'}} K_{y'})
}
\]
Now if $y \leq_{-} x \leq x'$, let $y_x \in Y_x$ be the element with $y
\leq y_x$. Then $y \leq y_x \leq_{-} x \leq x'$ and commutativity
follows from the previous two cases:
\[
\rho^{+}_{yx'} = \rho^{+}_{y_x x'} \rho^{+}_{y y_x} = \rho^{+}_{xx'}
\rho^{+}_{y_x x} \rho^{+}_{y y_x} = \rho^{+}_{xx'} \rho^{+}_{yx}
\]

The proof for $\rho^{-}$ in the cases $x' \leq x \leq_{+} y_x$ and $x'
\leq x \leq_{+} y$ is similar. Here we also use fact that $\vphi_{x'x}$
is an isomorphism to pick $y_{x'} = \vphi_{x'x}^{-1}(y_x)$ as the
unique element $y_{x'} \in Y_{x'}$ with $y_{x'} \leq y_x$.

\subsubsection{Construction of the natural transformations
$R^{+} R^{-} \to [1] \to R^{-} R^{+}$}
Observe that
\begin{align*}
&(\pmb{\xi}^{+} \pmb{\xi}^{-})_y = \xi_y[1]
&& (\pmb{\xi}^{-} \pmb{\xi}^{+})_y = \xi_y[1] \\
&(\pmb{\xi}^{+} \pmb{\xi}^{-})_x = \xi_{121} \circ \xi_{Y_x, x} &&
(\pmb{\xi}^{-} \pmb{\xi}^{+})_x = \xi_{212} \circ \xi_{x, Y_x}
\end{align*}
where $\xi_{121}$ and $\xi_{212}$ are the formulas defined
in~\eqref{e:D121}.

Let $\pmb{\nu}$ be the formula inducing the translation and define
$\eps^{+-} : \pmb{\xi}^{+} \pmb{\xi}^{-} \to \pmb{\nu}$, $\eps^{-+} :
\pmb{\nu} \to \pmb{\xi}^{-} \pmb{\xi}^{+}$ by
\begin{align*}
&\eps^{+-}_y : \xi_y[1] \xrightarrow{\left( \begin{smallmatrix} 1
\end{smallmatrix}\right)} \xi_y[1] \\
&\eps^{+-}_x : \xi_{121} \circ \xi_{Y_x,x} \xrightarrow{\beta_2 \circ
\xi_{Y_x,x}} \xi_2[1] \circ \xi_{Y_x,x} = \xi_x[1] \\
&\eps^{-+}_y : \xi_y[1] \xrightarrow{\left( \begin{smallmatrix} 1
\end{smallmatrix}\right)} \xi_y[1] \\
&\eps^{-+}_x : \xi_x[1] = \xi_1 \circ \xi_{x, Y_x}
\xrightarrow{\alpha_1 \circ \xi_{x,Y_x}} \xi_{212} \circ \xi_{x,Y_x,x}
\end{align*}
where $\xi_1$ and $\xi_2$ are as in~\eqref{e:xi12} and $\alpha_1$ and
$\beta_2$ are as in Proposition~\ref{p:Rpm1}. The proof of that
proposition also shows that $\eps^{+-}$ and $\eps^{-+}$ are morphisms
of formulas and induce natural transformations between functors, which
are quasi-isomorphisms.

\subsection{Proof of Corollary~\ref{c:X1Z}}
\label{ssec:X1Z}

Let $X$ and $Z$ be posets, and let $Y = \one \oplus Z$. Denote by $1
\in Y$ the unique minimal element and consider the map $f: X \to Y$
defined by $f(x)=1$ for all $x \in X$. Then
\begin{align*}
(X \sqcup Y, \leq^f_+) \simeq X \oplus \one \oplus Z && (X \sqcup Y,
\leq^f_-) \simeq \one \oplus (X + Z)
\end{align*}
hence by Corollary~\ref{c:XYorderf}, $X \oplus \one \oplus Z$ and $\one
\oplus (X + Z)$ are universally derived equivalent.


\begin{thebibliography}{1}

\bibitem{BGP73}
{\sc Bern{\v{s}}te{\u\i}n, I.~N., Gel{\cprime}fand, I.~M., and
Ponomarev,
  V.~A.}
\newblock Coxeter functors, and {G}abriel's theorem.
\newblock {\em Uspehi Mat. Nauk 28}, 2(170) (1973), 19--33.

\bibitem{Happel88}
{\sc Happel, D.}
\newblock {\em Triangulated categories in the representation theory of
  finite-dimensional algebras}, vol.~119 of {\em London Mathematical Society
  Lecture Note Series}.
\newblock Cambridge University Press, Cambridge, 1988.

\bibitem{Ladkani06}
{\sc Ladkani, S.}
\newblock On derived equivalences of categories of sheaves over finite posets.
\newblock {\tt arXiv:math.RT/0610685}.

\bibitem{Rickard89}
{\sc Rickard, J.}
\newblock Morita theory for derived categories.
\newblock {\em J. London Math. Soc. (2) 39}, 3 (1989), 436--456.

\bibitem{Rickard91}
{\sc Rickard, J.}
\newblock Derived equivalences as derived functors.
\newblock {\em J. London Math. Soc. (2) 43}, 1 (1991), 37--48.

\end{thebibliography}

\def\cprime{$'$}

\end{document}